\documentclass[aos,MSNbibl,seceqn,dvips]{arximspdf}
\usepackage{dcolumn}
\usepackage{graphicx}

\doi{10.1214/13-AOS1094} 
\volume{41}
\issue{2}
\pubyear{2013}
\firstpage{670}
\lastpage{692}

\makeatletter
\newcolumntype{d}[1]{D{.}{.}{#1}}
\def\cal{\mathcal}
\newcommand{\ka}{\kappa}

\newtheorem{lemma}{Lemma}[section]
\newproclaim{remark}{Remark}
\newproclaim{remarks}{Remarks}
\newtheorem{theorem}{Theorem}
\makeatother

\setattribute{abstract}{width}{290pt}

\begin{document}
\begin{frontmatter}

\title{Sequential multi-sensor change-point detection\thanksref{T1}}
\runtitle{Sequential multi-sensor change-point detection}
\thankstext{T1}{Supported in part by NSF Grant 1043204 and Stanford General Yao-Wu Wang Graduate \mbox{Fellowship}.}

\begin{aug}
\author[A]{\fnms{Yao} \snm{Xie}\corref{}\ead[label=e1]{yao.c.xie@gmail.com}}
\and
\author[B]{\fnms{David} \snm{Siegmund}\ead[label=e3]{siegmund@stanford.edu}}
\runauthor{Y. Xie and D. Siegmund}
\affiliation{Duke University and Stanford University}
\address[A]{Department of Electrical and Computer Engineering\\
Duke University\\
Durham, North Carolina, 27705\\
USA\\
\printead{e1}}

\address[B]{Department of Statistics\\
Stanford University\\
Stanford, California, 94305\\
USA\\
\printead{e3}}
\end{aug}

\received{\smonth{6} \syear{2012}}
\revised{\smonth{11} \syear{2012}}


\begin{abstract}
We develop a mixture procedure to monitor parallel streams of data for
a change-point that affects only a subset of them,
without assuming a spatial structure relating the
data streams to one another. Observations
are assumed initially to be independent standard normal random variables.
After a change-point the observations in a subset of the
streams of data have nonzero mean values.
The subset and
the post-change means are unknown.
The procedure we study uses stream specific generalized likelihood ratio
statistics, which are combined
to form an overall detection statistic
in a mixture model that hypothesizes an assumed
fraction $p_0$ of affected data streams.
An analytic expression is obtained for the average run length (ARL)
when there is no change and is shown by simulations to be very
accurate. Similarly, an approximation for the expected detection
delay (EDD) after a change-point is also obtained.
Numerical examples are given to compare the suggested procedure
to other procedures for unstructured problems and in one case
where the problem is assumed to have a well-defined geometric
structure. Finally we discuss sensitivity of the procedure
to the assumed value of $p_0$ and suggest a generalization.\vspace*{-2pt}
\end{abstract}

%

\end{frontmatter}

\section{Introduction}

Single sequence problems of change-point detection have a long history
in industrial quality control, where an observed process is assumed
initially to
be in control and at a change-point becomes out of control. It is
desired to detect the change-point with as little delay as possible,
subject to the constraint that false detections occurring before the
true change-point are very rare. Outstanding early
contributions are due to Page~\cite{Page1954,Page1955},
Shiryaev~\cite{Shiryaev1963} and Lorden
\cite{Lorden1971}.

We assume there
are parallel streams of data subject to change-points. More precisely
suppose that for each $n = 1, \ldots, N$, we make observations
$y_{n,t}$, $t = 1,2,\ldots.$ The observations\vadjust{\goodbreak} are mutually independent within
and across data streams.
At a certain time $\ka$,
there are
changes in the distributions of observations made at
a subset $\mathcal{N} \subset\{1, \ldots, N\}$ of cardinality
$|\mathcal{N}|
\leq N$.
Also denote by $\mathcal{N}^c$ the set of
unaffected data streams. The change-point $\ka$, the subset
$\mathcal{N}$ and its size,
and the size of the changes are all unknown.
As in the case of a single sequence, $N = 1$, the
goal is to detect the change-point as soon as possible after it occurs,
while keeping the frequency of false alarms as low as possible.
In the change-point detection literature, a surrogate for the
frequency of false alarms is
the average-run-length~(ARL), defined to be the expected time
before incorrectly announcing a change of distribution when none has occurred.

It may be convenient to imagine that the data streams represent
observations
at a collection of $N$ sensors and that the change-point is the onset
of a localized signal
that can be detected by sensors in the neighborhood of the signal.
In this paper we assume for the most part that the problem
is \textit{unstructured} in the sense that
we do not assume a model that relates the changes seen at the
different sensors.
An example of an unstructured problem is the model for
anomaly detection in computer networks developed in
\cite{Levy-LeducRoueff2009}.
For other discussions of
unstructured problems and applications, see
\cite{TartakovskyVeeravalli2008,Mei2009,Chen2010,PetrovRozovskiiTarkakovsky2003}.

At the other extreme
are \textit{structured} problems where there exists a profile
determining the relative
magnitudes of the changes observed by different sensors, say, according to
their distance from the location of a signal
(e.g.,~\cite{SiegmundYakir2008b,ShafieSigalSiegmund2003}).
A problem that is potentially
structured, may behave more like an unstructured problem if the
number of sensors is small and/or they are irregularly placed,
so their distances from one another are large compared to
the point spread function of the signals. Alternatively,
local signals may be collected at the relatively widely
dispersed hubs of small, star-shaped subnetworks,
then condensed and transmitted to a central processor, thus in effect
removing the
local structure.

The detection problem of particular interest in this paper
involves the case that $N$ is large and $|\mathcal{N}|$
is relatively small.
To achieve efficient detection, the detection procedure should use
insofar as possible only information from affected sensors and
suppress noise from the unaffected sensors.

In analogy to the well-known CUSUM statistic (e.g., Page~\cite{Page1954,Page1955},
Lorden~\cite{Lorden1971}),
Mei~\cite{Mei2009} recently proposed a multi-sensor
procedure based on sums of the
CUSUM statistic from individual sensors. He then compares the sum with a
suitable threshold to determine a stopping rule.
While the distributions of the data, both before
and after the change-point, are completely
general, they are also assumed to be completely known. The method is
shown to minimize asymptotically the expected detection delay (EDD) for
a given
false alarm rate, when the threshold value (and hence the
constraint imposed by the ARL) becomes infinitely large.
The procedure fails to be asymptotically optimal when the
specified distributions are incorrect.
Tartakovsky and Veeravalli\vadjust{\goodbreak} proposed a different procedure
\cite{TartakovskyVeeravalli2008} that sums the local likelihood ratio
statistic before forming CUSUM statistics. They also assume the
post-change distributions are completely prescribed. Moreover, both
procedures assume the
change-point is observed by all sensors. When only a subset of sensors
observe the change-point, these procedures include noise from the
unaffected sensors
in the detection statistic, which may lead to long detection delays.

In this paper, we develop a \textit{mixture}
procedure that achieves good detection performance in the case of an
unknown subset of affected sensors and incompletely specified
post-change distributions. The key feature of the proposed procedure is
that it incorporates an assumption about the fraction of affected sensors
when computing the detection statistic. We assume
that the individual observations are independent and normally
distributed with unit variance, and that the changes occur in their
mean values. At the $t$th vector of observations
$(y_{n,t}, n = 1, \ldots, N)$, the mixture procedure first computes a
generalized likelihood ratio (GLR) statistic for each individual sensor
under the assumption that a change-point has occurred at $k \leq t$.
The local GLR statistics are combined via
a mixture model that has the effect of soft thresholding the
local statistics according to an hypothesized fraction of
affected sensors, $p_0$.
The resulting local statistics are summed and compared with a detection
threshold. To characterize the performance of our proposed procedure, we
derive analytic approximations for its ARL and EDD,
which are evaluated by comparing the approximations to simulations.
Since simulation of the ARL is quite time consuming, the
analytic approximation to the ARL proves very useful in determining a
suitable detection threshold. The proposed procedure is then compared
numerically to competing procedures and is shown to be very competitive.
It is also shown to be reasonably robust
to the choice of $p_0$, and methods are suggested to increase
the robustness to mis-specification of $p_0$.

Although we assume throughout that the observations
are normally distributed, the model can be generalized to an
exponential family of distributions
satisfying some additional regularity conditions.

The remainder of the paper is organized as follows. In Section
\ref{Sec:model} we establish our notation and formulate the problem
more precisely.
In Section~\ref{Sec:Proc} we review several detection procedures and
introduce the proposed mixture procedure. In Section~\ref{Sec:ARL} we derive
approximations to the ARL and EDD of the mixture procedure,
and we demonstrate the accuracy of these approximations numerically.
Section~\ref{Sec:NumEg} demonstrates by
numerical examples that
the mixture procedure performs well compared to other procedures in the
unstructured problem. In Section~\ref{parallel} we suggest a ``parallel'' procedure
to increase robustness regarding the hypothesized fraction of affected
data streams, $p_0$. In Section~\ref{sec7} we also compare the mixture procedure
to that suggested in~\cite{SiegmundYakir2008b} for a structured
problem,
under the assumption that the assumed structure is correct.
Finally Section~\ref{Sec:conclusions} concludes the paper with some
discussion.\vadjust{\goodbreak}

\section{Assumptions and formulation}\label{Sec:model}

Given $N$ sensors, for each $n = 1, 2, \ldots,\break  N$, the observations from
the $n$th sensor are given by $y_{n,t}$, $t = 1, 2, \ldots.$ Assume that
different observations are mutually independent and normally distributed
with unit variances. Under the hypothesis of no change, they have
zero means. Probability and expectation in this case are denoted by
$\mathbb{P}^\infty$ and
$\mathbb{E}^\infty$, respectively. Alternatively,
there exists a change-point $\kappa$, $0 \leq\ka< \infty$, and a subset
$\mathcal{N}$ of $\{1, 2, \ldots, N\}$, having cardinality $|\mathcal
{N}|$, of
observations
\textit{affected} by the change-point.
For each $n \in\mathcal{N}$, the observations $y_{n,t}$
have means equal to $\mu_n > 0$ for all
$t > \ka$, while observations from the unaffected sensors
keep the same standard normal distribution. The probability and
expectation in this case are denoted by $\mathbb{P}^\ka$ and
$\mathbb{E}^\ka$, respectively. In particular,
$\kappa= 0$ denotes an immediate change.
Note that probabilities and expectations depend
on $\mathcal{N}$ and the values of $\mu_n$, although this
dependence is suppressed in the notation. The fraction of affected sensors
is given by $p = |\mathcal{N}|/N$.

Our goal is to define a stopping rule $T$ such that for all
sufficiently large prescribed
constants $c > 0$,
$\mathbb{E}^{\infty}\{T\}\geq c$, while asymptotically
$\mathbb{E}^\ka\{T- \ka|T > \ka\}$ is a minimum. Ideally, the
minimization would hold uniformly in the various unknown
parameters: $\ka$, $\mathcal{N}$ and the $\mu_n$. Since this is
clearly impossible,
in Section~\ref{Sec:NumEg} we will compare different procedures
through numerical examples computed under various hypothetical
conditions.\vspace*{-3pt}

\section{Detection procedures}\label{Sec:Proc}

Since the observations are independent, for an assumed value of the
change-point $\ka= k$ and sensor $n \in\mathcal{N}$, the log-likelihood
of observations
accumulated by time $t>k$ is given by
%
\begin{equation}
\label{LR} \ell_n (t, k,\mu_n) =\sum
_{i = k+1}^t \bigl(\mu_n y_{n,i} -
\mu_n^2/2 \bigr).
\end{equation}

We assume that each
sensor is affected by the change with probability
$p_0$ (independently from one sensor to the next).
The global log likelihood of all $N$ sensors is
%
\begin{equation}
\sum_{n=1}^N \log \bigl(1-p_0
+ p_0 \exp\bigl[\ell_n (t, k,\mu _n)\bigr]
\bigr) \label{mixloglik}.
\end{equation}
Expression (\ref{mixloglik}) suggests several change-point detection rules.

One possibility is to set $\mu_n$ equal to a nominal change, say
$\delta> 0$, which would
be important to detect, and define the stopping rule
%
\begin{equation}
\label{mixdelta} T_1 = \inf \Biggl\{t\dvtx \max
_{0\leq k \leq t} \sum_{n=1}^N \log
\bigl(1-p_0 + p_0 \exp\bigl[\ell^+_n (t, k,
\delta)\bigr]\bigr) \geq b \Biggr\},
\end{equation}
where $x^+$ denotes the positive part of $x$.
Here thresholding by the positive part plays the role of dimension
reduction by limiting the current considerations only to
sequences that appear to be affected by the change-point.\vadjust{\goodbreak}

Another possibility is to replace $\mu_n$ by its maximum likelihood
estimator, as follows.
The maximum likelihood estimate of the post-change mean as a function
of the current number of
observations $t$ and putative change-point location $k$ is given by
%
\begin{equation}
\hat{\mu}_n = \Biggl(\sum_{i=k+1}^t
y_{n,i} \Biggr)^+\Big/(t-k).
\end{equation}
Substitution into (\ref{LR}) gives the log generalized likelihood
ratio (GLR) statistic.
Putting
%
\begin{eqnarray}
\label{S_U} S_{n,t} & =& \sum_{i=1}^t
y_{n,i},
\nonumber
\\[-8pt]
\\[-8pt]
\nonumber
U_{n,k,t} & = &(t-k)^{-1/2}(S_{n,t} - S_{n,k}),
\end{eqnarray}
we can write the log GLR as
%
\begin{equation}
\ell_n(t,k,{\hat\mu_n}) = \bigl(U_{n,k,t}^+
\bigr)^2/2. \label{GLR_stats}
\end{equation}
We define the stopping rule
%
\begin{equation}
T_{2} = \inf \Biggl\{t\dvtx \max_{0\leq k < t}\sum
_{n=1}^N \log \bigl(1-p_0+
p_0 \exp\bigl[\bigl(U_{n,k,t}^+\bigr)^2/2\bigr]
\bigr) \geq b \Biggr\}.\label{Tmix}
\end{equation}

\begin{remark*} In what follows we use a \textit{window limited} version
of (\ref{Tmix}), where the maximum is restricted to
$m_0 \leq t - k < m_1$ for suitable $m_0 < m_1$. The role of~$m_1$ is two-fold. On the one hand it reduces the memory requirements
to implement the stopping rule, and on the other it effectively
establishes a minimum level of change that we want to detect.
For asymptotic theory given below, we assume that $b \rightarrow\infty$,
with $m_1/b$ also diverging. More
specific guidelines in selecting $m_1$ are discussed in~\cite{Lai1996window}.
In the numerical examples that follow, we take $m_0 = 1.$
In practice a slightly larger value can be used to provide
protection against outliers in the data, although it may delay
detection in cases involving very large changes.
\end{remark*}

The detection rule (\ref{Tmix}) is motivated by the suggestion of
\cite{ZhangYakirSiegmund2010}
for a similar fixed sample change-point detection problem.

For the special case $p_0 = 1$, (\ref{Tmix}) becomes the (global)
GLR procedure,
which for $N = 1$ was studied by~\cite{SiegmundVenkatraman1995}.
It is expected to be efficient if the change-point affects a large
fraction of the sensors.
At the other extreme, if only one or a very small number of
sensors is affected by the change-point,
a reasonable procedure would be
%
\begin{equation}
T_{\mathrm{max}} = \inf \Bigl\{t\dvtx \max_{0 \leq k < t }\max
_{1\leq n \leq
N} \bigl(U_{n,k,t}^+\bigr)^2/2 \geq b
\Bigr\}.\label{max_proc}
\end{equation}
The stopping rule $T_{\mathrm{max}}$ can also be
window limited.\vadjust{\goodbreak}

Still other possibilities are suggested by the observation that
a function of $y$
of the form $ \log[1-p_0 + p_0 \exp(y)]$ is large only if $y$ is large,
and then this function is approximately equal to $[y + \log(p_0)]^+$.
This suggests the stopping rules
%
\begin{equation}
\label{Tdelta} T_3 = \inf \Biggl\{t\dvtx \max_{0\leq k < t}
\sum_{n=1}^N \bigl[\ell
_n(t,k,\delta) + \log(p_0)\bigr]^+ \geq b \Biggr\}
\end{equation}
and
%
\begin{equation}
T_4 = \inf \Biggl\{ t\dvtx \max_{0\leq k < t} \sum
_{n=1}^N \bigl[\bigl(U^+_{n,k,t}
\bigr)^2/2 + \log(p_0)\bigr]^+ \geq b \Biggr\},
\end{equation}
or a suitably window limited version.

Mei~\cite{Mei2009} suggests the stopping rule
%
\begin{equation}
T_{\mathrm{Mei}} = \inf \Biggl\{t\dvtx \sum_{n=1}^N
\max_{0 \leq k < t} \ell_n(t,k,\delta) \geq b \Biggr\},
\label{mei}
\end{equation}
which simply adds the classical CUSUM statistics for the different
sensors. Note that this procedure does not involve the
assumption that all distributions affected by the change-point
change simultaneously. As we shall see below, this
has a negative impact on the efficiency of the procedure in our
formulation, although
it might prove beneficial in differently formulated problems.
For example, there may be a time delay before
the signal is perceived at different sensors, or there may be
different signals occurring at different times in the proximity of
different sensors.
In these problems, Mei's procedure, which allows
changes to occur at different times, could be useful.

The procedure suggested by Tartakovsky and Veeravalli
\cite{TartakovskyVeeravalli2008} is defined by the stopping rule
%
\begin{equation}
T_{\mathrm{TV}} \triangleq\inf \Biggl\{t\dvtx \max_{0 \leq k < t}\sum
_{n=1}^N \ell_n(t,k,\delta)
\geq b \Biggr\}.
\end{equation}
This stopping rule resembles $T_3(p_0)$ with $p_0 = 1$, but with one
important difference.
After a change-point the statistics of the
unaffected sensors have negative drifts that tend to cancel the positive
drifts from the affected sensors. This can lead to a large
EDD. Use of the positive part, $[\ell_n(t,k,\delta)]^+$, in the
definitions of our stopping rules is designed to avoid this problem.

Different thresholds $b$ are required for each of these
detection procedures to meet the ARL requirement.

\section{Properties of the detection procedures}\label{Sec:ARL}

In this section we develop theoretical properties of the
detection procedures $T_1$ to $T_4$, with emphasis on $T_2$ and
the closely related $T_4$.
We use two standard performance metrics:\vadjust{\goodbreak}
(i) the expected value of the stopping time when there is
no change, the average run length or ARL;
(ii) the expected detection delay (EDD), defined to be
the expected stopping time in the extreme case where
a change occurs immediately at $\ka= 0$.
The EDD provides an upper bound on
the expected delay after a change-point until detection occurs when the change
occurs later in the sequence
of observations.
The approximation to the ARL will be shown below to be very accurate,
which is fortunate since its simulation can be quite time
consuming, especially for large $N$. Accuracy of our approximation for
the EDD is variable, but fortunately this parameter is usually easily
simulated.

\subsection{Average run length when there is no change}

The ARL is the expected value of the stopping time $T$ when there
is no change-point. It will be convenient to use the following
notation. Let
$g(x)$ denote a twice continuously differentiable increasing function
that is bounded below at $- \infty$ and grows sub-exponentially at
$ + \infty$. In what
follows we consider explicitly $g(u)$. With an additional argument
discussed below the results also apply to $g(u^+)$.
Let
%
\begin{equation}
\psi(\theta) = \log\mathbb{E} \bigl\{\exp\bigl[\theta g(U)\bigr]\bigr\},
\end{equation}
where $U$ has
a standard normal distribution. Also let
%
\begin{equation}
\gamma(\theta) = \tfrac{1}{2}\theta^2 \mathbb{E}\bigl\{\bigl[
\dot{g}(U)\bigr]^2 \exp\bigl[\theta g(U) - \psi(\theta)\bigr]\bigr\},
\end{equation}
where the dot denotes differentiation.
Let
%
\begin{equation}
\label{abbrev} H(N,\theta) = \frac{\theta
[2\pi\ddot{\psi}(\theta)]^{1/2}}{ \gamma(\theta) N^{1/2}} \exp\bigl\{N\bigl[\theta\dot{
\psi}(\theta) - \psi(\theta)\bigr]\bigr\}.
\end{equation}
Denote the standard normal density function by $\phi(x)$ and its distribution
function by $\Phi(x)$.
Also let
$\nu(x) = 2x^{-2}\exp [-2\sum_{1}^\infty n^{-1}\Phi
(-|x|n^{1/2}/2 ) ]$; cf.~\cite{Siegmund1985}, page 82.
For numerical purposes a simple, accurate approximation is given by
(cf.~\cite{SiegmundYakir2007})
\[
\nu(x) \approx\frac{(2/x)[\Phi(x/2) - 0.5]}{(x/2)\Phi(x/2) + \phi(x/2)}.
\]

\begin{theorem}\label{th1}
Assume that
$N \rightarrow\infty$ and $b \rightarrow\infty$ with
$b/N$ fixed. Let $\theta$ be defined by
$\dot{\psi}(\theta) = b/N$.
For the
window limited stopping rule
%
\begin{equation}
T = \inf \Biggl\{t\dvtx \max_{0\leq k < t}\sum
_{n=1}^N g(U_{n,k,t}) \geq b \Biggr\}
\label{genstoprule}
\end{equation}
with $m_1 = o(b^r)$ for some
positive integer $r$, we have

\begin{equation}
\mathbb{E}^\infty\{T\} \sim H(N, \theta)\Big/ \int_{[2N\gamma(\theta
)/m_1]^{1/2}}^{[2N\gamma(\theta)/m_0]^{1/2}}
y \nu^2(y) \,dy. \label{ET}
\end{equation}
\end{theorem}

\begin{remark*}
The integrand in the approximation is integrable at
both $0$ and $\infty$ by virtue of the
relations $\nu(y) \rightarrow1$ as $y \rightarrow0$, and
$\nu(y) \sim2/y^2$ as $y \rightarrow\infty.$
\end{remark*}

The following calculations illustrate the essential features of
approximation~(\ref{ET}). For detailed proofs
in similar problems, see~\cite{SiegmundVenkatraman1995} (where additional
complications arise because
the stopping rule there is not window limited)
or~\cite{SiegmundYakir2008b}.
From arguments similar to those used in
\cite{ZhangYakirSiegmund2010}, we
can show that
%
\begin{eqnarray}
\label{SL} &&\mathbb{P}^\infty\{T \leq m\}
\nonumber\\
&&\qquad=\mathbb{P}^\infty \Biggl\{ \max_{t\leq m, m_0\leq t - k \leq m_1} \sum
_{n=1}^N g(U_{n,k,t}) \geq b \Biggr\}
\nonumber
\\[-8pt]
\\[-8pt]
\nonumber
&&\qquad\sim N^2 e^{-N[\theta\dot{\psi}(\theta) - \psi(\theta)]} \bigl[2\pi N \ddot{\psi}(\theta)
\bigr]^{-1/2}|\theta|^{-1} \gamma^2(\theta) \\
&&\hspace*{31pt}{}\times \int
_{m_0/m}^{m_1/m} \nu^2 \bigl(\bigl[2N \gamma(
\theta )/(mt)\bigr]^{1/2} \bigr) (1-t)\,dt/t^2.\nonumber
\end{eqnarray}
Here it is assumed that $m$
is large, but small enough that the right-hand side
of~(\ref{SL}) converges to 0 when $b\rightarrow\infty$.
Changing variables in the integrand and using
the notation (\ref{abbrev}), we can re-write this
approximation as
%
\begin{equation}
\mathbb{P}^\infty\{T \leq m\} \sim m \int_{[2N\gamma(\theta)/m_1]^{1/2}}
^{[2N\gamma(\theta)/m_0]^{1/2}} y \nu^2(y)\,dy/H(N,\theta). \label{LD}
\end{equation}
From the arguments in~\cite{SiegmundVenkatraman1995} or
\cite{SiegmundYakir2008b} (see also~\cite{Aldous1988}),
we see that
$T$ is
asymptotically exponentially
distributed and is uniformly integrable. Hence if $\lambda$ denotes
the factor
multiplying $m$ on the right-hand side
of (\ref{LD}), then for still larger $m$, in the range
where $m \lambda$ is bounded
away from 0 and $\infty,$
$\mathbb{P}^\infty\{T \leq m\} - [1 - \exp(- \lambda m)] \rightarrow0.$
Consequently
$\mathbb{E}^\infty\{T\} \sim\lambda^{-1}$,
which is equivalent to (\ref{ET}).

\begin{remarks*} (i) The result we have used from
\cite{ZhangYakirSiegmund2010} was motivated by a problem involving
changes that could be
positive, or negative, or both; and in that paper it was assumed that the
function $g(u)$ is
twice continuously differentiable.
The required smoothness is not satisfied by the composite
functions of principal
interest here, of the form $g(u^+)$.
However, (i) the required smoothness is required only in the
derivation of some of the constant factors, not the exponentially small
factor,
and (ii) the second derivative that appears in the
derivation in~\cite{ZhangYakirSiegmund2010} can be eliminated from the
final approximation by an integration by parts.
As a consequence, we can approximate the indicator of $u > 0$ by $\Phi
(ru)$ and
use in place of $u^+$
the smooth function $\int_{- \infty}^u \Phi(rv) \,dv
= u \Phi(ru) + r^{-1} \phi(ru)$, which converges uniformly to
$u^+$ as $r \rightarrow\infty.$
Letting $r \rightarrow\infty$\vadjust{\goodbreak}
and interchanging limits produce~(\ref{SL}). An alternative approach would be simply to define
$g(u)$ to be appropriate for a one-sided change while having the
required smoothness in $u$.
An example is
$g(u) = \log[1-p_0 + p_0 \exp(u^2\Phi(ru)/2)]$, which sidesteps
the technical issue, but seems less easily motivated.

(ii)
The fact that all the stopping times studied in this paper
are asymptotically exponentially
distributed when there is no change
can be very useful.
(A simulation illustrating this property in the case of $T_2$ is
given in Section~\ref{sec4.3}.)
To simulate the ARL, it is not
necessary to simulate the process until the stopping time $T$, which
can be
computationally time consuming, but only until a time $m$ when
we are able to estimate $\mathbb{P}^\infty\{T \leq m\}$ with a
reasonably small percentage error. For
the numerical examples given later,
we have occasionally used this shortcut with the value of $m$ that makes
this probability
0.1 or 0.05.

(iii) Although the mathematical assumptions involve large values of $N$,
some numerical experimentation for $T_2(p_0)$ shows that
(\ref{ET}) gives roughly
the correct values even
for $N = 1$ or 2. For $p_0 = 1$ (\ref{ET}) provides numerical
results similar to those given for
the generalized likelihood ratio
statistic in~\cite{SiegmundVenkatraman1995}.

(iv) Theorem~\ref{th1} allows us to approximate the ARL for
$T_2$ and $T_4$. The stopping rule $T_{\mathrm{max}}$ is
straightforward to handle, since the minimum of $N$ independent
exponentially distributed random variables is itself exponentially
distributed. The stopping rules $T_1$ and $T_3$, where $g$
is composed with $\ell_{t,k,\delta}^+$, can be handled by
a similar argument with one important difference. Now
the cumulant generating function $\psi(\theta)$ depends on
$w = t-k$, so the equation defining $\theta$ must
be solved for each value of $w$, and the resulting approximation
summed over possible values of $w$. Fortunately only a few terms
make a substantial contribution to the sum, except when $\delta$
is very small. For the results reported below, the additional
amount of computation is negligible.
\end{remarks*}
\subsection{Expected detection delay}\label{Sec:DD}

After a change-point occurs, we are interested in the
expected number of additional
observations required for detection.
For the detection rules considered in this paper, the maximum
expected detection
delay over $\ka\geq0$ is attained at $\ka= 0$. Hence we consider this
case.

Here we are unable to consider stopping times defined by
a general function~$g$, so we consider the specific functions involved
in the definitions of $T_2$ and $T_4$.
Let $g(u, p_0) = \log(1-p_0 + p_0 \exp[(u^+)^2/2])$ or
$[(u^+)^2/2 + \log(p_0)]^+$, and
let $U$ denote a standard normal random variable.
Recall that ${\cal N}$ denotes the set of sensors at which
there is a change, $|\mathcal{N}|$ is the cardinality of this set
and $p = |\mathcal{N}|/N$ is the true fraction of sensors that are
affected by the change.
For each $n \in{\cal N}$ the mean value changes from 0
to $\mu_n > 0$, and for $n \in\mathcal{N}^c$ the distribution remains
the same as before the change-point. Let
%
\begin{equation}
\Delta= \biggl(\sum_{n \in{\cal N}} \mu_n^2
\biggr)^{1/2}.
\end{equation}

Note that the Kullback--Leibler divergence of a vector of
observations after the change-point from a vector of observations before
the change-point
is $\Delta^2/2$, which determines the asymptotic rate of growth of
the detection statistic after the change-point. Using Wald's
identity~\cite{Siegmund1985}, we see to a first-order approximation that
the expected detection delay is $2b/\Delta^2$, provided that the
maximum window
size, $m_1$, is large compared to this quantity. In the following derivation
we assume $m_1 \gg2 b/ \Delta^2$.

In addition, let
%
\begin{equation}
\tilde{S}_t \triangleq\sum_{i=1}^t
z_i \label{process}
\end{equation}
be a random walk where the increments $z_i$ are independent and
identically distributed with mean $\Delta^2/2$ and variance $\Delta^2$.
Let $\tau= \min \{t\dvtx \tilde{S}_t > 0 \}$.
Our approximation to the expected detection delay given
below depends on two related quantities. The first is
%
\begin{equation}
\rho(\Delta) = \tfrac{1}{2} \mathbb{E}\bigl\{{\tilde S}_\tau^2
\bigr\}/\mathbb {E}\{\tilde{S}_\tau\}
\end{equation}
for which exact computational expressions
and useful approximations are
available in~\cite{Siegmund1985}.
In particular,
%
\begin{equation}
\rho(\Delta) = \mathbb{E}\bigl\{z_1^2\bigr\}/\bigl(2
\mathbb{E}\{z_1\}\bigr) - \sum_{i=1}^\infty
i^{-1} \mathbb{E}\bigl\{\tilde{S}_i^-\bigr\} =
\Delta^2/4 +1- \sum_{i=1}^\infty
i^{-1} \mathbb{E}\bigl\{\tilde{S}_i^-\bigr\},\hspace*{-35pt}
\end{equation}
where $(x)^- = -\min\{x, 0\}$. The second quantity is $\mathbb{E}\{
\min_{t\geq0}\tilde{S}_t\}$, which
according to (Problem 8.14 in~\cite{Siegmund1985}) is given by
%
\begin{equation}
\mathbb{E} \Bigl\{\min_{t\geq0}\tilde{S}_t \Bigr\} =
\rho(\Delta ) - 1 - \Delta^2/4.
\end{equation}

The following approximation refines the first-order result
for the expected detection delay. Recall that $\mathbb{E}^0$
denotes expectation when the change-point $\kappa= 0.$

\begin{theorem}\label{th2}
Suppose $b \rightarrow\infty$, with other parameters held fixed.
Then for $T = T_2$ or $T_4$,
%
\begin{eqnarray}
\label{DD} \mathbb{E}^0\{T\}
&=& 2\Delta^{-2} \Bigl[ b+ \rho(\Delta) - |\mathcal{N}|\log
p_0 - |\mathcal{N}|/2 + \mathbb{E} \Bigl\{\min_{t \geq0}{
\tilde S}_t \Bigr\}
\nonumber
\\[-8pt]
\\[-8pt]
\nonumber
&&\hspace*{81pt}{} -\bigl(N - |\mathcal{N}|\bigr)\mathbb{E} \bigl\{g(U,
p_0)\bigr\} + o(1) \Bigr].
\end{eqnarray}
\end{theorem}

The following calculation provides the ingredients for
a proof of (\ref{DD}). For details in similar problems
involving a single sequence, see~\cite{PollakSiegmund1975}
and~\cite{SiegmundVenkatraman1995}. For convenience
we assume that $T = T_{2}$, but there is almost no difference
in the calculations when $T= T_4$. Let $k_0 = b^{1/2}$.
For $k < T- k_0$, we can write the detection statistic at the stopping
time $T$ as follows, up\vadjust{\goodbreak} to a term that tends to zero exponentially fast
in probability:
%
\begin{eqnarray}
\label{decomp} Z_{k, T}&=& \sum_{n=1}^N
g(U_{n,k,T}, p_0)
\nonumber\\[-2pt]
&=& \sum_{n\in\mathcal{N}} g(U_{n,k,T}, p_0) +
\sum_{n\in\mathcal
{N}^c} g(U_{n,k,T}, p_0)
\nonumber\\[-2pt]
&=& \sum_{n\in\mathcal{N}} \log \biggl(p_0\exp \bigl
\{\bigl(U_{n,k,T}^+\bigr)^2/2 \bigr\} \biggl[1+
\frac{1-p_0}{p_0}\exp \bigl\{-\bigl(U_{n,k,T}^+\bigr)^2/2 \bigr
\} \biggr] \biggr)\nonumber\\[-2pt]
&&{} + \sum_{n\in\mathcal{N}^c} g(U_{n,k,T},
p_0)
\nonumber\\[-2pt]
&=& \sum_{n\in\mathcal{N}} \bigl[\log p_0 +
\bigl(U_{n,k,T}^+\bigr)^2/2 \bigr] + \sum
_{n\in\mathcal{N}^c} g(U_{n,k,T}, p_0) \\[-2pt]
&&{}+ \sum
_{n\in
\mathcal{N}} \log \biggl(1+\frac{1-p_0}{p_0}\exp \bigl\{ -
\bigl(U_{n,k,T}^+\bigr)^2/2 \bigr\} \biggr)
\nonumber\\[-2pt]
&=& |\mathcal{N}| \log p_0 + \sum_{n\in\mathcal{N}}
\bigl(U_{n,k,T}^+\bigr)^2/2 + \sum
_{n\in\mathcal{N}^c} g(U_{n,k,T}, p_0) + o(1)
\nonumber\\[-2pt]
&=& |\mathcal{N}|\log p_0 + \sum_{n\in\mathcal{N}}
\bigl[(S_{n,T}-S_{n,k})^+\bigr]^2/2(T-k)\nonumber \\[-2pt]
&&{}+ \sum
_{n\in\mathcal{N}^c} g(U_{n,k,T}, p_0) + o(1).\nonumber
\end{eqnarray}
The residual term $\sum_{n\in\mathcal{N}} \log (1+(1-p_0)\exp
\{-(U_{n,k,T}^+)^2/2 \}/p_0 )$ tends to zero
exponentially fast when $b\rightarrow
\infty$ because when $b\rightarrow
\infty$, $T\rightarrow b/\Delta$, and $n\in\mathcal{N}$,
$(U_{n,k,T}^+)^2$ grows on the order of $\mu_n^2 (T-k) > \mu_n^2 k_0
= \mu_n^2 \sqrt{b}$.\vspace*{1pt}

We then use the following simple identity to decompose the second term
in~(\ref{decomp}) for the affected sensors into two parts:
%
\begin{eqnarray}
\label{key_id} \bigl(S_{n,t}^+\bigr)^2/2t &=&
S_{n,t}^2/2t - \bigl(S_{n,t}^-\bigr)^2/2t
\nonumber
\\[-8pt]
\\[-8pt]
\nonumber
&=& \mu_n(S_{n,t} - \mu_n t/2) +
(S_{n,t} - \mu_n t)^2/2t -
\bigl(S_{n,t}^-\bigr)^2/2t.
\end{eqnarray}
From the preceding discussion, we see that $\max_{0\leq k < T - k_0}
Z_{k,T}$ is on the order of $b$, while $\max_{T-k_0 \leq k < T} Z_{k,
T}$ is on the order of $k_0 = b^{1/2}$.
Hence with overwhelming probability the max over all $k$ is attained for
$k < T-k_0$, so from (\ref{key_id}) and (\ref{decomp}) we have
%
\begin{eqnarray}
\label{19}
&&\max_{0 \leq k < T} Z_{k, t}
\nonumber\\
&&\qquad=\max_{0 \leq k < T-k_0} \sum_{n=1}^N
g(U_{n,k,T}, p_0) + o(1)
\nonumber\\
&&\qquad= |\mathcal{N}|\log p_0\nonumber\\
&&\qquad\quad{}+ \max_{0 \leq k < T-k_0} \biggl[\sum
_{n\in
\mathcal
{N}} \mu_n \bigl[(S_{n,T}
-S_{n,k}) - (T-k)\mu_n/2\bigr]
\nonumber\\
&&\hspace*{94pt}{}+ \sum_{n\in\mathcal{N}} \bigl[(S_{n,T} -
S_{n,k})-(T-k)\mu _n\bigr]^2/\bigl[2(T-k)\bigr]
\nonumber\\
&&\hspace*{94pt}{} - \bigl[(S_{n,T}-S_{n,k})^-\bigr]^2/2(T-k)
\nonumber
\\[-8pt]
\\[-8pt]
\nonumber
&&\hspace*{215pt}{}+ \sum_{n\in\mathcal{N}^c} g(U_{n,k,T}, p_0)
\biggr] + o(1)
\\
&&\qquad= |\mathcal{N}|\log p_0+\sum
_{n\in\mathcal{N}} \mu_n (S_{n,T} - T\mu
_n/2 )
\nonumber\\
&&\qquad\quad{}+\max_{0 \leq k < T-k_0} \biggl[ - \sum_{n\in\mathcal{N}}
\mu _n (S_{n,k} - k\mu_n/2 ) \nonumber\\
&&\hspace*{92pt}{}+ \sum
_{n\in\mathcal{N}} \bigl[(S_{n,T} - S_{n,k})-(T-k)
\mu_n\bigr]^2/\bigl[2(T-k)\bigr]
\nonumber\\
&&\hspace*{59pt}\qquad\quad{}      - \sum_{n\in\mathcal{N}} \bigl[(S_{n,T}-S_{n,k})^-
\bigr]^2/\bigl[2(T-k)\bigr] \nonumber\\
&&\hspace*{213pt}{}+ \sum_{n\in\mathcal{N}^c} g
(U_{n,k,T}, p_0 ) \biggr] + o(1).\nonumber
\end{eqnarray}

The following lemma forms the basis for the rest of the
derivation. The proof is omitted here;
for details see~\cite{YaoXie2011} (or~\cite{SiegmundVenkatraman1995}
for the special case $N = 1$).

\begin{lemma}\label{lemma1}
For $k_0 = b^{1/2}$, asymptotically as $b\rightarrow\infty$
\begin{eqnarray*}
&&\max_{0 \leq k < T} \biggl[- \sum_{n\in\mathcal{N}}
\mu_n (S_{n,k} - k\mu_n/2 ) %
+ \sum
_{n\in\mathcal{N}} \frac{[(S_{n,T} - S_{n,k})-(T-k)\mu_n]^2}{2(T-k)} %
 \\
 &&\hspace*{102pt}{}- \sum_{n\in\mathcal{N}} \frac{[(S_{n,T} - S_{n,k})^-]^2}{2(T-k)} %
+
\sum_{n\in\mathcal{N}^c} g(U_{n,k,T}, p_0)
\biggr]
\\
&&\qquad= \sum_{n\in\mathcal{N}} (S_{n,T} -T
\mu_n)^2/2T + \sum_{n\in\mathcal{N}^c}
g(U_{n,0,T}, p_0)\\
&&\qquad\quad{} + \max_{0 \leq k < k_0} \biggl[- \sum
_{n\in\mathcal{N}} \mu _n(S_{n,k} - k
\mu_n/2) \biggr] + o_p(1),
\end{eqnarray*}
where $o_p(1)$ converges to 0 in probability.
\end{lemma}

By taking expectations in (\ref{19}), letting $b\rightarrow\infty$ and
using Lemma~\ref{lemma1}, we have
%
\begin{eqnarray}
\label{main} &&  \mathbb{E}^0 \Biggl\{\max_{0 \leq k < T}\sum
_{n=1}^N g(U_{n,k,T},
p_0) \Biggr\}
\nonumber\\
&&\qquad=  \mathbb{E}^0 \biggl\{|\mathcal{N}|\log p_0 + \sum
_{n\in\mathcal
{N}} \mu _n(S_{n,T} - T
\mu_n/2) + \sum_{n\in\mathcal{N}} \frac{(S_{n,T} -
T\mu_n)^2}{2T}
\nonumber
\\[-8pt]
\\[-8pt]
\nonumber
&&\hspace*{55pt}{} +
\sum_{n\in\mathcal{N}^c} g(U_{n,0,T}, p_0)
+ \max_{0 \leq k < k_0} \biggl[-\sum
_{n\in\mathcal{N}} \mu_n(S_{n,k} -k
\mu_n/2) \biggr] \biggr\} \\
&&\qquad\quad{}+ o(1).\nonumber
\end{eqnarray}
We will compute each term on the right-hand side of (\ref{main})
separately. We will need the lemma due to Anscombe and Doeblin (see
Theorem 2.40 in~\cite{Siegmund1985}), which states that the standardized randomly
stopped sum
of random variables are asymptotically normally distributed under quite
general conditions.

\begin{longlist}[(iii)]
\item[(i)]
By Wald's identity~\cite{Siegmund1985},
%
\begin{equation}
\mathbb{E}^0 \biggl\{\sum_{n\in\mathcal{N}}
\mu_n(S_{n,T} - T\mu _n/2) \biggr\} =
\mathbb{E}^0\{T\}\Delta^2/2.
\end{equation}
\item[(ii)] By the Anscombe--Doeblin lemma, $(S_{n,T} - T\mu
_n)/T^{1/2}$ is asymptotically
normally distributed with zero mean and unit variance. Hence\break
$\sum_{n\in\mathcal{N}} (S_{n,T} - T\mu_n)^2/T$ is asymptotically a
sum of independent $\chi^2_1$ random variables, so
%
\begin{equation}
\mathbb{E}^0 \biggl\{\sum_{n\in\mathcal{N}}
(S_{n,T} - T\mu _n)^2/2T \biggr\} = |
\mathcal{N}|/2 + o(1).
\end{equation}
\item[(iii)] Similarly,
%
\begin{equation}
\mathbb{E}^0 \biggl\{\sum_{n\in\mathcal{N}^c}
g(U_{n,0,T}, p_0) \biggr\} \rightarrow\bigl(N - |\mathcal{N}|\bigr)
\mathbb{E}^0\bigl\{g(U, p_0)\bigr\}.
\end{equation}
\item[(iv)] The term
$ -\sum_{n\in\mathcal{N}} \mu_n(S_{n, k} - \mu_n k/2)$ ($k \geq
0$) is a
random walk with~negative drift $- \Delta^2/2$ and variance $\Delta^2$. Hence
$\mathbb{E}^0  \{
\max_{0 \leq k< k_0} [-\sum_{n\in\mathcal{N}}
\mu_n(S_{n, k} -k\mu_n/2)] \}$
converges to the expected minimum of this
random walk, which has the same distribution as
$\min_{t \geq0} {\tilde S}_t$ defined above.
\end{longlist}

Having evaluated the right-hand side of (\ref{main}), we
now consider the left-hand side, to which we will apply a
nonlinear renewal theorem. This requires that we write the process of
interest
as a random walk and a relatively slowly varying remainder, and
follows standard lines by using
a Taylor series approximation to show that for large values
of $t$ and bounded values of $k$ (cf.~\cite{SiegmundVenkatraman1995,PollakSiegmund1975}, and the argument already given above)
the asymptotic growth of $\sum_{n=1}^N g(U_{n, k,t}, p_0)$
for $t > \kappa$ is governed by the
random walk $\sum_{n\in\mathcal{N}}\mu_n(S_{n,t} - t\mu_n/2)$,
which has
mean value $t \Delta^2/2$ and variance $t \Delta^2$.
By writing
%
\begin{equation}
\mathbb{E}^0 \Biggl\{\max_{0 \leq k < T} \sum
_{n=1}^N g(U_{n,k,T}, p_0) \Biggr
\} = b + \mathbb{E}^0 \Biggl\{\max_{0 \leq k < T} \sum
_{n=1}^N g(U_{n,k,T},
p_0) - b \Biggr\},\hspace*{-35pt}
\end{equation}
and using nonlinear renewal theory to evaluate the expected overshoot
of the
process of (\ref{process}) over the boundary
(\cite{Siegmund1985}, Chapter IX),
we obtain
%
\begin{equation}
\mathbb{E}^0 \Biggl\{\max_{0 \leq k < T} \sum
_{n=1}^N g(U_{n,k,T}, p_0) - b
\Biggr\} \rightarrow\rho(\Delta).
\end{equation}

\begin{remarks*}
(i) Although the proof of Theorem~\ref{th2} follows the
pattern of
arguments given previously in the case $N = 1$, unlike that case
where the asymptotic approximation is surprisingly accurate even
when the EDD is relatively small, here
the accuracy is quite variable.
The key element in the derivation is the asymptotic linearization of
$g(U^+_{n,k,t}, p_0)$
for each $n \in{\cal N}$ into a term involving a random walk
and a remainder. A simple test for conditions when the approximation
will be reasonably accurate is to compare the exact value of
$\mathbb{E}^0\{Z_{0,t}\}$, which is easily evaluated by numerical
integration, to the expectation of the linearized
approximation, then take $t$ large enough to make these
two expectations approximately equal. If such a value of $t$ makes
the expectations less than or equal to $b$, the approximation of the theorem
will be reasonably accurate. Indeed the preceding argument is
simply an elaboration of these equalities at $t = T$ combined with
Wald's identity to extract $\mathbb{E}^0\{T\}$ from the random walk,
and numerous technical steps to approximate the nonnegligible
terms in the remainders. For a crude, but quite reliable approximation
that has no mathematically precise justification that we can see,
choose $t$ to satisfy $\mathbb{E}^0\{Z_{0,t}\} = b$. Fortunately
the EDD is easily simulated when it is small, which is where
problems with the analytic approximation arise.

(ii) In principle the same method can be used to approximate the expected
detection delay of $T_1$ or $T_3$. In some places the analysis is substantially
simpler, but in one important respect it is more complicated. In
the preceding argument, for $n \in{\cal N}^c$ the term involving
the expected value of $g(U_{n,0,T},p_0)$ is very simple, since
$U_{n,0,T}^2$ has asymptotically a $\chi^2$ distribution.
For the stopping rules $T_1$ and $T_3$, the
term $g(\ell_{n,0,T}, p_0)$ does not have a limiting distribution,
and in fact for $n \in{\cal N}^c$ it converges to
0 as $b \rightarrow\infty.$ However, there
are typically a large number of these terms, and in many
cases $T$ is relatively
small, nowhere near its asymptotic limit. Hence it would be unwise
simply to replace this expectation by 0. To a crude first-order approximation
$T \sim b/[\delta_0(\sum_{n \in{\cal N}} \mu_n - \delta_0/2)]$ =
$t_0$, say.
Although it is not correct mathematically speaking, an often reasonable
approximation can be obtained by using the term
$-(N-|\mathcal{N}|) \mathbb{E}^0\{\ell_{n,0,t_0}\}$
to account for the statistics associated with sequences
unaffected by the change-point.
Some examples are included in the numerical examples in Table~\ref{table:methods_comparison}.
\end{remarks*}
\subsection{Accuracy of the approximations}\label{sec4.3}

We start with examining the accuracy
of our approximations for the ARL and the EDD in
(\ref{ET}) and (\ref{DD}).
For a Monte Carlo experiment we use $N = 100$ sensors, $m_1 = 200$ and
$\mu_n = 1$ for all affected data streams.
The comparisons for different
values of $p_0$ between the theoretical and Monte Carlo ARLs obtained
from 500 Monte Carlo trials are given in
Tables~\ref{table:ARL2} and~\ref{table:ARL4}, which show that the
approximation
in (\ref{ET})
is quite accurate.

\begin{table}
\tablewidth=250pt
\caption{ARL of $T_2(p_0), m_1 = 200$}\label{table:ARL2}
\begin{tabular*}{250pt}{@{\extracolsep{\fill}}lccc@{}}
\hline
$\bolds{p_0}$& $\bolds{b}$ & \textbf{Theory} & \textbf{Monte Carlo}
\\
\hline
0.3 & 31.2 & \phantom{0.}5001 & \phantom{0.}5504 
\\
0.3 & 32.3 & 10,002 & 10,221 
\\[3pt]
0.1& 19.5 & \phantom{0.}5000 & \phantom{0.}4968 
\\
0.1 & 20.4 & 10,001 & 10,093 
\\[3pt]
0.03 & 12.7 & \phantom{0.}5001 & \phantom{0.}4830
\\
0.03 & 13.5 & 10,001 & \phantom{0.}9948 
\\
\hline
\end{tabular*}
\end{table}

\begin{table}[b]
\tablewidth=250pt
\caption{ARL of $T_4(p_0)$, $m_1 = 200$}\label{table:ARL4}
\begin{tabular*}{250pt}{@{\extracolsep{\fill}}lccc@{}}
\hline
$\bolds{p_0}$& $\bolds{b}$ & \textbf{Theory} & \textbf{Monte Carlo}
\\
\hline
0.3 & 24.0 & 5000 & 5514
\\
0.1& 15.1 & 5000 & 5062 
\\
0.03 & 10.8 & 5000 & 5600
\\
\hline
\end{tabular*}
\end{table}

\begin{figure}

\includegraphics{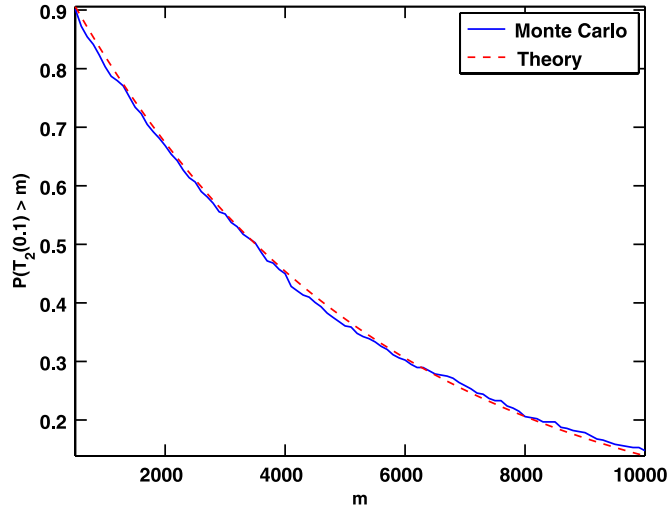}

\caption{The tail probability $\mathbb{P}\{T_2(0.1) > m\}$.
Approximate theoretical values are obtained from~(\protect\ref{ET});
numerical values are
obtained from 500 Monte Carlo trials.}
\label{fig:tail}
\end{figure}

Figure~\ref{fig:tail} illustrates the fact that $T_2(0.1)$ is
approximately exponentially
distributed.

Results for the EDD
obtained from 500 Monte Carlo trials are given in
Table~\ref{table:DD}. Although the approximation seems reasonable,
it does not appear to
be as accurate as the approximation for the ARL. Since
the EDD
requires considerably less computational effort to simulate
and needs to be known only roughly when we choose design parameters
for a particular problem,
there is less value to an accurate analytic approximation.

\begin{table}[b]
\tabcolsep=3pt
\caption{EDDs of $T_2(p_0)$ and $T_4(p_0)$ with ARL
$\approx5000$, $\mu= 1$, and $m_1 = 200$}\label{table:DD}
\begin{tabular*}{\textwidth}{@{\extracolsep{\fill}}ld{1.2}d{2.1}d{2.1}d{2.1}d{2.1}@{}}
\hline
\multicolumn{1}{@{}l}{$\bolds{p}$} & \multicolumn{1}{c}{$\bolds{p_0}$} &
\multicolumn{1}{c}{\textbf{Theory} $\bolds{T_2(p_0)}$} & \multicolumn{1}{c}{\textbf{Monte Carlo} $\bolds{T_2(p_0)}$} & \multicolumn{1}{c}{\textbf{Theory}
$\bolds{T_4(p_0)}$} &\multicolumn{1}{c@{}}{\textbf{Monte Carlo} $\bolds{T_4(p_0)}$} \\
\hline
0.3 & 0.3 & 3.5 & 3.2 & 4.2 & 3.5\\
0.1 & 0.3 & 6.2 & 6.5 & 7.1 & 6.6 \\
0.3 & 0.1 & 5.2 & 3.6 & 5.1 & 4.1\\
0.1 & 0.1 & 7.2 & 6.7 &7.0 & 7.1 \\
0.03 & 0.1 & 13.9 & 14.3& 13.5 & 14.3 \\
0.03 & 0.03 & 13.9 & 14.2 & 13.7 & 14.6 \\
\hline
\end{tabular*}
\end{table}

%
\begin{table}
\tablewidth=250pt
\caption{Thresholds for ARL $\approx$ 5000, $m_1 = 200$}\label{ARL_compare_2}
\begin{tabular*}{250pt}{@{\extracolsep{\fill}}lcc@{}}
\hline
\textbf{Procedure} & $\bolds{b}$ & \textbf{Monte Carlo ARL} \\
\hline
Max & 12.8 & 5041\\
$T_2(1)$ & 53.5 & 4978 \\
$T_2(0.1)$ & 19.5 & 5000\\
Mei & 88.5 & 4997\\
$T_3(0.1,1)$ & 12.4 & 4948 \\
$T_3(1,1)$ & 41.6 & 4993 \\
\hline
\end{tabular*}
\end{table}

\begin{table}
\caption{EDD with $N = 100$ obtained from 500 Monte Carlo trials.
Thresholds for ARL 5000 are listed in Table \protect\ref{ARL_compare_2}.
Theoretical approximations for EDD are in parentheses}\label{table:methods_comparison}
\begin{tabular*}{\textwidth}{@{\extracolsep{\fill}}lcccc@{}}
\hline
\multicolumn{1}{@{}l}{$\bolds{p}$}&
\multicolumn{1}{c}{\textbf{Method}} &
\multicolumn{1}{c}{\textbf{EDD,} $\bolds{\mu=1}$} &
\multicolumn{1}{c}{\textbf{EDD,} $\bolds{\mu= 0.7}$} &
\multicolumn{1}{c@{}}{\textbf{EDD,} $\bolds{\mu= 1.3}$}\\
\hline
0.01 & max & 25.5 & 49.6 & 16.3 \\
& $T_2(1)$ & 52.3 (56.9) & 105.5 (114.6) & 32.9 (34.1) \\
& $T_2 (0.1)$ & 31.6 (32.5) & 59.4 (64.9) & 20.3 (19.7)\\
& Mei & 53.2 & 103.8 & 38.1 \\
&$T_3(0.1,1)$ & 29.1 (29.3) & 63.3 (59.0) & 19.1 (19.1) \\
& $T_3(1,1)$ & 82.0 (83.6) & 213.7 (193.5) & 53.3 (53.5) \\[3pt]
0.03 & max & 18.1 & 33.3 & 11.6 \\
& $T_2(1)$ & 18.7 (19.3) & 35.8 (38.4) & 12.6 (11.7) \\
& $T_2(0.1)$ & 14.2 (13.9) & 26.7 (27.5) & 9.3 (8.5) \\
& Mei & 23.0 & 41.6 & 16.4\\
& $T_3(0.1)$ & 13.4 & 26.9 & 9.2 \\
& $T_3(1)$ & 27.2 & 66.0 & 16.3 \\[3pt]
0.05 & max & 15.5 & 28.4 & 9.7 \\
& $T_2(1)$ & 12.2 (11.6) & 21.8 (23.0) & 7.9 (7.1)\\
& $T_2(0.1)$ & 10.4 (10.1) & 18.9 (19.9) & 6.9 (6.2) \\
& Mei & 15.7 & 26.9 & 11.4 \\
& $T_3(0.1,1)$ & 9.8 (9.8) &18.6 (21.4) & 7.0 (6.8) \\
& $T_3(1,1)$ & 15.5 (16.2) & 38.8 (39.8) & 9.0 (9.7)\\[3pt]
0.1 & max & 12.6 & 23.0 & 8.4 \\
& $T_2(1)$ & 6.7 (5.9) & 11.8 (11.3) & 4.7 (3.7) \\
& $T_2(0.1)$ & 6.7 (7.2) & 11.6 (14.1) & 4.6 (4.5) \\
& Mei & 9.6 & 15.4 & 7.4 \\
& $T_3(0.1,1)$ & 7.1 (7.6) & 11.9 (16.7) & 5.3 (5.3) \\
& $T_3(1,1)$ & 6.8 (7.3) & 15.7 (19.6) & 4.6 (4.5) \\[3pt]
0.3& max & 9.6 & 16.7 & 6.6\\
& $T_2(1)$ & 3.0 (2.0) & 4.4 (3.5) & 2.4 (1.4) \\
& $T_2(0.1)$ & 3.5 (5.2) & 5.6 (10.1) & 2.7 (3.3) \\
& Mei & 4.9 & 7.0 & 4.0\\
& $T_3(0.1,1)$ & 4.6 & 6.7 & 3.9 \\
& $T_3(1,1)$ & 3.0 & 4.3 & 2.5 \\[3pt]
0.5 & max & 8.6& 14.4 & 5.8\\
& $T_2(1)$ & 2.3 & 3.0 & 2.0\\
& $T_2(0.1)$ & 2.8 & 4.0 & 2.1\\
& Mei & 3.8 & 5.0 & 3.0\\
& $T_3(0.1,1)$ &4.0 &5.4 & 3.3 \\
& $T_3(1,1)$ & 2.3 & 3.0 & 2.0 \\[3pt]
1 & max & 7.2 & 12.1 & 5.1\\
& $T_2(1)$ & 2.0 & 2.0 & 2.0 \\
& $T_2(0.1)$ & 2.0 & 2.6 & 2.0\\
& Mei & 3.0& 3.4 & 2.3 \\
& $T_3(0.1,1)$ & 3.4 & 4.3 & 3.0\\
& $T_3(1,1)$ & 2.0 & 2.1 & 2.0\\
\hline
\end{tabular*}
\end{table}

We have performed considerably more extensive simulations that yield results
consistent with the small experiments reported in Tables
\ref{table:ARL2},~\ref{table:ARL4} and
\ref{table:DD}.
Since the parameter $p_0$ defining $T_2$ must be chosen subjectively,
it is interesting
to observe that
Table~\ref{table:DD} suggests these procedures are
reasonably robust with respect to the choice of $p_0$,
and choosing $p_0$ somewhat
too large seems less costly than choosing $p_0$ too small. More extensive
calculations bear out this observation. We return to the problem of choosing
$p_0$ in Section~\ref{parallel}.

\section{Numerical comparisons}\label{Sec:NumEg}

In this section, we compare the
expected detection delays for several procedures when their ARLs are all
approximately 5000.
The thresholds are given in Table~\ref{ARL_compare_2}, where we assume
$N = 100$, and $m_1 = 200$
for those procedures for which a limited window size is appropriate.
Procedure~(\ref{Tmix}) is
denoted by $T_2(p_0)$.
For Mei's procedure we put $\delta_n = 1.$ The procedures in
(\ref{Tdelta}) are denoted by
$T_3(p_0,\delta).$ Recall that $T_2(1)$ uses the generalized
likelihood ratio
statistic and $T_3(1,\delta)$ is similar to the procedure proposed by
Tartakovsky and
Veeravalli~\cite{TartakovskyVeeravalli2008}, but we have inserted the positive part to avoid the
problems mentioned in Section~\ref{Sec:Proc}.
The expected detection delays are obtained from 500 Monte Carlo trials
and are listed in Table~\ref{table:methods_comparison}.
For some entries, values from our asymptotic approximation are given in
parentheses.

Note that the max procedure (\ref{max_proc}) has the smallest
detection delay when $p = 0.01$,
but it has the largest delay for $p$ greater than 0.1.
The procedures defined by
$T_2$ and by $T_3$ are comparable. Mei's procedure performs well
when $p$ is large, but
poorly when $p$ is small.

\section{Parallel mixture procedure}\label{parallel}

The procedures considered above depend on a parameter $p_0$, which
presumably should be chosen to be close to the unknown true fraction $p$.
While Table~\ref{table:methods_comparison} suggests that the value of
$p_0$ is fairly robust when $p_0$ does not
deviate much from the true $p$,
to achieve robustness over a wider range of the unknown parameter $p$,
we consider a parallel procedure that combines several procedures, each using
a different $p_0$ to monitor a different range of $p$ values. The thresholds
of these individual procedures will be chosen so that they have the
same ARL.
For example, we can use two different values of $p_0$, say a small $p_0
= p_1$
and a large $p_0 = p_2$, and then choose thresholds $b_1$ and $b_2$ to
obtain the same ARL for these two procedures.
The parallel procedure claims a detection if at least one of the component
procedures reaches its threshold, specifically
%
\begin{equation}
T_{\mathrm{parallel}}\triangleq\min\bigl\{T_2(p_1),
T_2(p_2)\bigr\}.
\end{equation}

To compare the performance of the parallel procedure with that
of a single $T_2$, we
consider a case with $N=400$ and $m_1= 200$. For the single mixture procedure
we use the intermediate value $p_0 = 0.10$
and threshold value $b = 44.7$, so
$\mathbb{P}^\infty\{T_2 \leq1000\} \approx0.10$ and hence the
ARL $\approx 10\mbox{,}000$.
For the parallel procedure we consider the values $p_1 = 0.02$
and $p_2 = 0.33$. For the
threshold values $b_1 = 21.2 $ and $b_2 = 87.7$, respectively, we have
$\mathbb{P}^\infty\{T_2(p_i) \leq1000 \} \approx0.05$, $i = 1, 2$.
By the Bonferroni inequality
$\mathbb{P}^\infty\{ \min[T_2(p_1),T_2(p_2)] \leq1000 \} \leq0.1,$ so
conservatively
$\mathbb{E}^\infty\{T_{\mathrm{parallel}}\} \geq10\mbox{,}000.$
Table~\ref{table:parallel} shows that the expected
detection delays of the parallel procedure are usually smaller
than those of the
single procedure, particularly for very small or very large $p$.
Presumably these differences are magnified in problems involving
larger values of $N$, which have the possibility of still smaller
values of $p$.

%
\begin{table}
\tablewidth=250pt
\caption{Comparison of EDD, parallel and simple procedures}\label{table:parallel}
\begin{tabular*}{250pt}{@{\extracolsep{\fill}}lccc@{}}
\hline
$\bolds{p}$ & $\bolds{\mu}$ & $\bolds{T_2(0.1)}$\textbf{, EDD} & \textbf{Parallel, EDD} \\
\hline
0.1 & 0.7 & \phantom{0}6.5 & \phantom{0}6.4 \\
0.005 & 1.0 & 27.1 & 22.9\\
0.005 & 0.7 & 54.5 & 45.8\\
0.25 & 0.3 & 12.0 & 10.5\\
0.4 & 0.2 & 14.4 & 12.3 \\
0.0025 & 1.5 & 23.3 & 17.8\\
\hline
\end{tabular*}
\end{table}

Simulations indicate that because of dependence between the two
statistics used to define the parallel procedure, the ARL is actually
somewhat larger than the Bonferroni approximation suggested. Since the
parallel procedure becomes increasingly attractive in larger problems,
which provide more room for improvement over a single choice of
$p_0$, but which are also increasingly difficult to simulate,
it would be interesting to develop a more accurate theoretical
approximation for the ARL.

An attractive alternative to the parallel procedure would be to use
a weighted linear combination for different values of $p_0$
of the statistics used to define
$T_2$ or $T_3$. Our approximation for the ARL can be
easily adapted, but some modest numerical exploration suggests that
the expected detection delay is not improved as much as for the
parallel procedure.

\section{Profile-based procedure for structured problems}\label{sec7}

Up to now we have assumed there is no spatial structure relating
the change-point amplitudes
at difference sensors. In this section we will consider briefly
a \textit{structured problem},
where there is a parameterized profile of the amplitude of the signal
seen at each
sensor that is based on the distance of the sensor to the
source of the signal.
Assuming we have some knowledge about such a profile, we can incorporate
this knowledge into the definition of an appropriate detection
statistic. Our developments follow closely the analysis in \cite
{SiegmundYakir2008b}.

Assume the location of the $n$th sensor is given by its
coordinates $x_n$, $n = 1, \ldots, N$ at points in Euclidean
space, which for simplicity we take to be on an equi-spaced grid.
We assume that
the source is located in a region $\mathcal{D}$, which is a subset of
the ambient Euclidean space. In our example below we consider two
dimensional space, but three dimensions would also be quite reasonable.
Assume the change-point amplitude at the $n$th sensor is determined by
the expression
%
\begin{equation}
\mu_n = \sum_{m=1}^M
r_m \alpha_{z_m}(x_n), \label{prof_fun_def}
\end{equation}
where $M$ is the number of sources, $z_m \in\mathcal{D}$ is
the (unknown) spatial location of the $m$th source,
$\alpha_z(x)$ is the
profile function, and
the scalar $r_m$ is
an unknown parameter that measures the strength of the $m$th signal.
The profile function
describes how the signal strength of the $m$th point source
has decayed at the $n$th sensor. We assume some knowledge about this
profile function is available. For example, $\alpha_z(x)$ is often
taken to be a decreasing
function of the Euclidean distance between $z$ and~$x$.
The profile may also depend on finitely many parameters, such as the
rate of
decay of the function. See~\cite{Rabinowitz1994} or
\cite{ShafieSigalSiegmund2003} for examples in a fixed sample context.

If the parameters $r_m$ are multiplied by a positive
constant and the profile $\alpha_{z_m}(x_n)$ divided by the same
constant, the
values of $\mu_n$ do not change. To avoid this lack of identifiability,
it is convenient to assume that for all $z$ the
profiles have been standardized to have unit Euclidean norm, that is,
$\sum_{x} \alpha^2_z(x) = 1$ for all z.

\subsection{Profile-based procedure}
Under the assumption that there is at most one source,
say at $z$, for observations up to time $t$ with a change-point
assumed to equal~$k$,
the log likelihood function for observations from all sensors (\ref
{LR}) is
%
\begin{equation}
\ell(t, k, r, z) = \sum_{n=1}^N \bigl[r
\alpha_z(x_n) (S_{n,t} - S_{n,k}) -
r^2 (t-k) \alpha^2_z(x_n)/2
\bigr]. \label{likelihood_profile}
\end{equation}
When maximized with respect to $r$ this becomes
%
\begin{equation}
\frac{1}{2} \biggl[\biggl\{\sum_n
\alpha_z(x_n)U_{n,k, t}\biggr\}^+
\biggr]^2. \label{l_profile}
\end{equation}
Maximizing the function (\ref{l_profile}) with respect to the
putative change-point $k$ and the source location $z$, we obtain the
log GLR statistic and a profile-based
stopping rule of the form
%
\begin{equation}
T_{\mathrm{profile}} = \inf \biggl\{t\dvtx \max_{0 \leq k < t }\max
_{z \in\mathcal{D}} \biggl[\biggl\{\sum_n
\alpha_z(x_n) U_{n, k, t}\biggr\}^+
\biggr]^2\geq b \biggr\}.\label{profile_2}
\end{equation}
If the model is correct,
(\ref{profile_2}) is a matched-filter type of statistic.

\subsection{Theoretical ARL of profile-based procedure}\label{sec:ARL_profile}

Using the result presented in~\cite{SiegmundYakir2008b}, we can derive an approximation for the
ARL of the profile-based procedure. We consider
in detail a special case\vadjust{\goodbreak} where $d = 2$ and the profile is given
by a Gaussian function
%
\begin{equation}
\alpha_{z}(x) = \frac{1}{\sqrt{2\pi\beta}}e^{-(1/(4\beta))\Vert x-z\Vert ^2},\qquad x\in
\mathbb{R}^2, \beta> 0.\label{prof_def}
\end{equation}
The parameter
$\beta>0$ controls of rate of profile decay and is assumed known.
With minor modifications one could also maximize with respect to a
range of values of~$\beta$.

Although the sensors have been assumed to be located on the
integer lattice of two-dimensional Euclidean space, it will be
convenient as a very rough approximation to assume that summation over
sensor locations $x$
can be approximated by integration over the entire Euclidean space.
With this approximation, $\sum_x \alpha^2_z(x)$, which we have assumed
equals 1 for all $z$, becomes $\int_{\mathbb{R}^2} \alpha^2_z(x) \,dx$,
which by (\ref{prof_def}) is readily seen to be identically 1.
The approximation is reasonable if $\beta$ is large, so the effective
distance between points of the grid is small, and the
space ${\cal D}$, assumed to contain the signal, is well within
the set of sensor locations (so edge effects can be ignored and the
integration extended over all of ${\mathbb{R}^2})$.

It will be convenient to use the notation
%
\begin{equation}
\label{innerprod} \langle f, g\rangle= \int_{\mathbb{R}^2} f(x) g(x)
\,dx.
\end{equation}
Let $\dot{\alpha}_z$ denote the gradient of
$\alpha_z$ with respect to $z$. Then according to
\cite{SiegmundYakir2008b},
%
\begin{eqnarray}
\label{prof_theo} &&\mathbb{P}^{\infty}\{T_{\mathrm{profile}}\leq m
\}
\nonumber
\\
&&\qquad\sim m \exp(-b/2) (b/4 \pi)^{3/2}2^{1/2} \\
&&\qquad\quad{}\times\int
_{(b/m_1)^{1/2}}^{(b/m_0)^{1/2}}u \nu^2(u) \,du \int
_{\cal D} \bigl|\operatorname{det} {\bigl(\bigl\langle\dot{\alpha}_{z},
\dot{\alpha}_{z}^{\top}\bigr\rangle\bigr)\bigr|^{1/2}} \,dz.\nonumber
\end{eqnarray}
To evaluate the last integral in (\ref{prof_theo}),
we see from (\ref{prof_def}) that $\dot{\alpha}_z$ satisfies
%
\begin{equation}
\dot{\alpha}_{z}(x) = \alpha_z(x) (x-z)/(2\beta).
\end{equation}
Hence by (\ref{innerprod})
$\langle\dot{\alpha}_{z}, \dot{\alpha}_{z}^{\top}\rangle$
is a $2\times2$ matrix of integrals, which can be easily evaluated,
and its
determinant equals $1/(16\beta^4)$.
Hence the last integral in
(\ref{prof_theo}) equals
$|{\cal D}|/(4\beta^2)$
where $|{\cal D}|$ denotes the area of ${\cal D}$.
Arguing as above from the asymptotic exponentiality of
$T_{\mathrm{profile}}$, we find that an asymptotic
approximation for the average run length is given by
%
\begin{eqnarray}
&&\mathbb{E}^\infty\{T_{\mathrm{profile}}\}
\nonumber
\\[-8pt]
\\[-8pt]
\nonumber
&&\qquad \sim 16\bigl(2\pi^3
\bigr)^{1/2} \beta^2 b^{-3/2} \exp(b/2)\Big/ \biggl[ \int
_{(b/m_1)^{1/2}}^{(b/m_0)^{1/2}} u \nu^2(u) \,du\cdot |{\cal D}|
\biggr].
\end{eqnarray}

\subsection{Numerical examples}

In this section we briefly compare the unstructured detection procedure
based on $T_2$ with the profile-based procedure in the special case that
the assumed profile is correct.\vadjust{\goodbreak}

Assume that the profile is given by the Gaussian function (\ref{prof_def})
with parameter $\beta= 1$ and both procedures are
window-truncated with $m_0 = 1$, $m_1 = 100$.
The number of sensors is $N = 625$ distributed over a $25 \times25$
square grid with center
at the origin. In this situation,
approximately $p = 0.016$ sensors are affected. In the specification of
$T_2$, we take $p_0 = 0.05$.

The thresholds are chosen so that the average run lengths when there
is no change-point are approximately 5000.
Using (\ref{prof_theo}), we obtain
$\mathbb{P}^\infty\{T_{\mathrm{profile}} \leq250\} = 0.050$ for $b = 29.5$.
From 500 Monte Carlo trials we obtained
the threshold 26.3,
so the theoretical approximation appears to be slightly conservative.

To deal with a failure to know the true
rate of decay of the signal with distance, we could maximize over
$\beta$, say, for $\beta\in[0.5,5]$.
A suitable version of
(\ref{prof_theo}) indicates the threshold would be 33.8.
This slight increase to the threshold suggests that failure to know the
appropriate
rate of decay of the signal with distance leads to a relatively moderate
loss of detection efficiency.

%
\begin{table}
\tablewidth=265pt
\caption{Comparison of EDD, profile-based and unstructured procedures}\label{table:DD_large}
\begin{tabular*}{265pt}{@{\extracolsep{\fill}}lccc@{}}
\hline
& $\bolds{b}$ & \textbf{EDD} $\bolds{r = 1}$ & \textbf{EDD,} $\bolds{r = 1.5}$\\
\hline
Profile-based procedure & 26.3 & 25.6 & 12.3 \\
Unstructured procedure & 39.7 & 78.3 & 35.8 \\
\hline
\end{tabular*}
\end{table}

For comparisons of the EDD, we used for the
profile-based procedure the threshold 26.3, given by simulation, while for
$T_2(0.05)$ we used the analytic approximation, which our studies
have shown to be very
accurate.
Table~\ref{table:DD_large} compares the expected detection delay of the
profile-based procedure with that of the mixture procedure.
As one would expect from the precise modeling assumptions, the
profile-based procedure is substantially more powerful.

In many cases there will be only a modest scientific basis for the
assumed profile,
especially in multidimensional problems. The distance between sensors
relative to
the decay rate of the signal is also an important consideration. It
would
be interesting to compare the structured and the unstructured
problems when the
assumed profile
differs moderately or substantially from the true profile,
perhaps in the number of sources of the signals, their shape,
the rate of decay, or the locations of the sensors.

\section{Discussion}\label{Sec:conclusions}

For an unstructured multi-sensor change-point problem we have suggested
and compared a number of sequential detection procedures.
We assume that the pre- and post-change samples are normally distributed
with known variance and that both the post-change mean and the set of
affected sensors are unknown. For performance analysis, we have derived
approximations for the average run length (ARL) and the expected detection
delay (EDD),\vadjust{\goodbreak} and have shown that these approximations have reasonable accuracy.
Our principal procedure depends on the assumption that a known
fraction of sensors are affected by the change-point.
We show numerically that the procedures are
fairly robust with respect to discrepancies between the actual and the
hypothesized fractions, and we suggest a parallel procedure based on two
or more hypothesized fractions to increase this robustness.

In a structured problem, we have shown that knowledge of the
correct structure can be implemented to achieve large improvements
in the EDD. Since the assumed structure is usually at best only
approximately correct, an interesting
open question is the extent to which failure to hypothesize the appropriate
structure compromises these improvements. One possible method to
achieve robustness against inadequacy of the structured model would be
a parallel version of structured and unstructured detection.

%



\printaddresses

\end{document}